\documentclass[11pt]{article}
\usepackage{amssymb,amsfonts,amsmath,latexsym,epsf,tikz,url}
\usepackage[usenames,dvipsnames]{pstricks}
\usepackage{pstricks-add}
\usepackage{epsfig}
\usepackage{pst-grad} % For gradients
\usepackage{pst-plot} % For axes
\usepackage[space]{grffile} % For spaces in paths
\usepackage{etoolbox} % For spaces in paths
\usepackage{float}
\usepackage{multicol}
\usepackage{soul}
\usepackage{tikz}
\usepackage[colorinlistoftodos]{todonotes}
\usepackage{pgfplots}
\usepackage{mathrsfs}
\usepackage[colorlinks]{hyperref}
\usetikzlibrary{arrows}
\makeatletter % For spaces in paths
\patchcmd\Gread@eps{\@inputcheck#1 }{\@inputcheck"#1"\relax}{}{}
\makeatother

\newtheorem{thm}{Theorem}[section]
\newtheorem{cor}[thm]{Corollary}
\newtheorem{lem}[thm]{Lemma}

\newtheorem{defn}[thm]{Definition}
\newtheorem{obs}[thm]{Observation}
\newtheorem{rem}[thm]{\bf{Remark}}

\numberwithin{equation}{section}

\newcommand{\proof}{\noindent{\bf Proof.\ }}
\newcommand{\qed}{\hfill $\square$\medskip}

\begin{document}

\def\nt{\noindent}

\title{Fair coalition in graphs}

\author{	
Saeid Alikhani\footnote{Corresponding Author}, Abbas Jafari, Maryam Safazadeh
}

%\date{\today}

\maketitle

\begin{center}
	
Department of Mathematical Sciences, Yazd University, 89195-741, Yazd, Iran\\

\medskip 
{\tt alikhani@yazd.ac.ir,~ abbasjafaryb91@gmail.com, ~ msafazadeh92@gmail.com,}

\end{center}

%%%%%%%%%%%%%%ABSTRACT%%%%%%%%%%%%%%%%%%%%%%%%%%%%%%%%%%%%%%%%%%%%%%%%%%%%%%%%%%%%
\begin{abstract}
 Let $G=(V,E)$ be a simple graph. A dominating set of $G$ is a subset $D\subseteq V$ such that every vertex not in $D$ is adjacent to at least one vertex in $D$.
 The cardinality of a smallest dominating set of $G$, denoted by $\gamma(G)$, is the domination number of $G$.  For $k \geq 1$, a $k$-fair dominating set ($kFD$-set) in $G$, is a dominating set $S$ such that $|N(v) \cap D|=k$ for every vertex $ v \in V\setminus D$. A fair dominating set in $G$ is a $kFD$-set for some integer $k\geq 1$. We consider $1FD$-sets and define a fair coalition in a graph $G$ as a pair of disjoint subsets $A_1, A_2 \subseteq A$ that satisfy the following conditions: (a) neither $A_1$ nor $A_2$ constitutes a $1$-fair  dominating set of $G$, and (b) $A_1\cup A_2$ constitutes a $1$-fair dominating set of $G$.	
 A fair coalition partition of a graph $G$ is a partition $\Upsilon = \{A_1,A_2,\ldots,A_k\}$ of its vertex set   
 such that every set $A_i$ of $\Upsilon$ is either
 a singleton $1$-fair dominating set of $G$, or is not a $1$-fair dominating set of $G$ but forms a fair coalition with another non-$1$-fair dominating set $A_j\in \Upsilon$. 
 We define the fair coalition number of $G$ as the maximum cardinality of a fair coalition partition of $G$, and we denote it by $\mathcal{C}_f(G)$.
 We initiate the study of the fair coalition in graphs and obtain $\mathcal{C}_f(G)$  for some specific graphs. 
\end{abstract}


\medskip
\noindent{\bf AMS Subj.\ Class.:} 05C25, 05C60.

\section{Introduction}

Let $G=(V,E)$ be a simple graph. 
A set $D\subseteq V$ is a  dominating set, if every vertex in $V\backslash D$ is adjacent to at least one vertex in $D$. 
The  domination number $\gamma(G)$ is the minimum cardinality of a dominating set in $G$. A vertex $v$ in graph $G$ of order $n$ is called a full vertex of $G$, if $\deg(v)=n-1$. 

Haynes et al. \cite{6} first defined the concept of a coalition in graphs as two non-dominating sets whose union is dominating, and subsequently introduced the coalition partition and the coalition number. Their work established initial bounds for the coalition number and determined it for paths and cycles. Later research, such as in \cite{9}, expanded on these bounds by considering minimum and maximum degrees. The study of coalition graphs, where adjacency represents coalition formation, was initiated in \cite{7}, showing that all graphs can be coalition graphs. This was further explored for specific graph families like trees, paths, and cycles in \cite{8}. 

The $c$-partition problem has also been investigated in specific graph types, such as trees by Bakhshesh et al. \cite{3} and cubic graphs by Alikhani, Golmohammadi and  Konstantinova \cite{1}. Alikhani et al. have also studied variations of coalition partitions, including total \cite{total}  and connected \cite{2} coalitions.  Jafari, Alikhani and Bakhshesh \cite{Abbas} for $k$-coalitions, Mojdeh et al. for perfect  and edge \cite{Mojdeh1,Mojdeh2} coalitions.

A domatic partition is a partition
of the vertex set into dominating sets. The maximum cardinality of a domatic partition is called the domatic number, denoted by $d(G)$. The domatic number of a graph was introduced in 1977 by  Cockayne and  Hedetniemi \cite{5}.

A dominating set $D$ in a graph $G$ is an $i$-fair dominating set (or $iFD$-set) if every vertex $v\in V\setminus D$ has exactly $i$ neighbors in $D$, for some integer $i\geq 1$. The $i$-fair domination number of $G$, denoted by $fd_i(G)$, is defined as the minimum cardinality of an $iFD$-set. An $iFD$-set that achieves this minimum cardinality is termed an $fd_i(G)$-set. More broadly, a fair dominating set (abbreviated FD-set) is any $iFD$-set for some $i\geq 1$. The fair domination number of a graph $G$ (if $G$ is not the empty graph), symbolized as $fd(G)$, corresponds to the minimum cardinality among all $FD$-sets. 
If $G$ is the empty graph on $n$ vertices, then $fd(G)$  is conventionally defined as $n$. From these definitions, it follows that for any graph $G$ of order $n$, $\gamma(G)\leq fd(G)\leq n$, and the equality $fd(G)=n$ holds precisely when $G=\overline{K_n}$. Caro, Hansberg, and Henning \cite{Henning}  have made notable contributions to this area, including demonstrating that for a disconnected graph $G$ (without isolated vertices) of order $n\geq 3$, $fd(G)\leq n-2$, and constructing families of graphs that achieve this bound. They further established that for a tree $T$ of order $n\geq 2$, $fd(T)\leq \frac{n}{2}$ , with equality if and only if $T$ is a specific type of tree,  $T=T'\circ K_1$ which is the corona product of a tree $T'$ and $K_1$. 
The enumerative aspects of $1$-fair dominating sets have been explored in research by Alikhani and Safazadeh  \cite{Fair1,Fair2}.
We denote the $1$-fair domination number of $G$ which is the minimum cardinality of  $1FD$-set  by $\gamma_f(G)$.

\medskip 
We consider $1$-fair dominating sets, introduce and initiate the study of the fair coalition in graphs and obtain $\mathcal{C}_f(G)$ for some specific graphs in Section 2. We obtain the fair coalition number of cubic graphs of order at most  $10$ in Section 3. Finally, we conclude
the paper in Section 4. 

\section{Introduction to fair coalition} 
We begin this section by formally defining fair coalitions within graph theory. Subsequently, we derive several bounds for this newly introduced parameter. The section concludes with a determination of the fair coalition number for both path graphs and cycle graphs.

We first  define a fair domatic and a fair coalition  and then we establish some results. 

\begin{defn}
	A $1$-fair dominating set ($1FD$-set) in $G$, is a dominating set $S$ such that $|N(v) \cap D|=1$ for every vertex $ v \in V\setminus D$.
	A fair domatic partition  is a partition
	of the vertex set into $1$-fair dominating sets. The maximum cardinality of a fair domatic partition  is called the fair domatic number, denoted by $d_f(G)$.
\end{defn}

\begin{defn}[Fair coalition]
	A fair coalition in a graph $G$ consists of two disjoint sets
	$A_1$ and $A_2$ of vertices of $G$, neither of which is a $1$-fair dominating set but whose union
	$A_1\cup A_2$ is a $1$-fair dominating set of $G$. 
\end{defn}
Let us introduce  fair coalition partition for a graph $G$.
\begin{defn}[Fair coalition partition]
	A fair coalition partition, abbreviated $fc$-partition, of a graph $G$ refers to a vertex partition $\Upsilon=\{A_1,\ldots, A_k\}$,    
	such that every set $A_i$ of $\Upsilon$ is either
	a singleton $1$-fair dominating set of $G$, or is not a $1$-fair dominating set of $G$ but forms a fair coalition with another non-$1$-fair dominating set $A_j\in \Upsilon$.  The fair coalition number of $G$, denoted by $\mathcal{C}_f(G)$, refers to the largest possible number of members in a $fc$-partition of $G$. A $fc$-partition of $G$ of cardinality $\mathcal{C}_f(G)$ is called a $\mathcal{C}_f(G)$-partition.
\end{defn}

\subsection{Some bounds}

First we establish a relation between the fair coalition number $C_f(G)$ and the fair domatic number $d_{f}(G)$ as follows.

\begin{thm}\label{codo}
	If $G$ is a graph of order $n\geq 3$ without full vertices, then $C_f(G) \ge 2d_{f}(G)$.
\end{thm} 
\begin{proof} 
	Let $G$ has a fair domatic partition $\mathcal{S}=\{S_1, S_2,\ldots, S_k\}$ with $d_{f}(G)=k$. Since $G$ has no vertices of degree $n-1$ then $|S_i|>1$ for any $i$.  Without loss of generality we assume that the sets $\{S_1,S_2,\ldots,
	S_{k-1}\}$ are minimal $1$-fair dominating sets. Indeed, if for some $i$, the set $S_i$ is not minimal, we find a subset $S'_i\subseteq S_i$ that is a minimal $1$-fair dominating set, and add the remaining vertices to the set $S_k$. Note that if we partition a minimal $1$-fair dominating set with more than one element into two non-empty sets, we obtain two non-$1$-fair dominating sets that together form a fair coalition. As a result, we divide each non singleton set $S_i$ into two sets $S_{i,1}$ and $S_{i,2}$ that form a fair coalition. This gives us a new partition $\mathcal{S}'$ consisting of non-fair dominating sets that pair with some other non-$1$-fair dominating set in $\mathcal{S}'$  form a fair coalition. 
	
	We now check the $1$-fair dominating set $S_k$. 
	
	If $S_k$ is a minimal $1$-fair dominating set, we divide it into two non-$1$-fair dominating sets, add these sets to $\mathcal{S}'$, and obtain a fair  coalition partition of order at least $2k$. Then, since $k=d_{f}(G)$, $C_f(G)\ge 2d_{f}(G)$.
	
	If $S_k$ is not a minimal $1$-fair dominating set, we aim to get a subset $S'_k\subseteq S_k$ that holds this condition. Again, we use the strategy on partitioning $S'_k$ into two non-$1$-fair dominating sets giving together a fair coalition. Afterwards, we define $S''_k$ as the complement of $S'_k$ in $S_k$, and append $S'_{k,1}$ and $S'_{k,2}$ to $\mathcal{S}'$. If $S''_k$ can merge
	with any non-$1$-fair dominating set to form a fair coalition, one can obtain a fair coalition partition of a cardinality at least $2k+1$ by adding $S''_K$ to $\mathcal{S}'$. Then, $C_f(G)\ge 2d_{f}(G)+1$. However, if $S''_k$ can not form a fair coalition with any set in $\mathcal{S}'$, we eliminate $S'_{k,2}$ from $\mathcal{S}'$ and add the set $S'_{k,2}\cup S''_k$ to $\mathcal{S}'$. This leads to a fair coalition partition of a cardinality at least $2k$. Then, $C_f(G)\ge 2d_{f}(G)$.
	
	Due to the above arguments, we  conclude that $SC(G)\ge 2d_{f}(G)$. \qed
\end{proof}

\begin{rem}
	There are examples which show that the bound in Theorem \ref{codo} is sharp. It suffices to consider $P_4$ with $V=\{v_1,v_2,v_3,v_4\}$ 
	and
	$E=\{v_1v_2,v_2v_3,v_3v_4\}$. It is easy to see that $C_f(P_{n})=4$. Consider two sets $S_1=\{v_1,v_4\}$ and $S_2=\{v_2,v_3\}$. These two sets $S_1$ and $S_2$ are $1$-fair dominating sets and we have found a partition of $V(P_4)$  into two disjoint fair dominating set $S_1$ and $S_2$, so $d_f(P_4)=2$. 
\end{rem} 

We continue the study of this parameter. For instance we obtain some bounds for $\mathcal{C}_f(G)$ based  the $1$-fair domination number, i.e., $\gamma_f(G)$.

\begin{thm}\label{upper}
	Let $G$ be a graph with order $n$ and $1$-fair domination number $\gamma_f$. Then
	$$\mathcal{C}_f(G)\leq n-\gamma_f+2.$$
\end{thm}
\begin{proof}
	For graphs $G$ with $\gamma_f=1$, the inequality is obviously true. If $\gamma_f(G)\geq 2$, then $G$ has no full vertex. Let $t= {\mathcal C}_f(G)$ and let $\Upsilon=\{A_1,A_2,\dots, A_t\}$
	be a $C_f(G)$-partition of $G$. So, we have
	\begin{eqnarray}
	n=|A_1|+|A_2|+\dots+|A_t|.
	\end{eqnarray}
	Without loss of generality, assume that $A_1$ and $A_2$ form a fair coalition. Then
	$|A_1|+|A_2|\geq \gamma_f(G)$. Combining this with (1), we obtain
	$$n\geq |A_1|+|A_2|+t-2.$$
	Therefore we have the results. \qed 
\end{proof} 

\begin{rem}
	The equality in Theorem \ref{upper}
	holds for many graphs. For example consider the complete bipartite graph 
	$K_{3,3}$. You can find a $1$-fair dominating set by selecting one vertex from each partition set. So $\gamma_f(K_{3,3})=2$. In the other hand $\mathcal{C}_f(K_{3,3})=6$ (see Theorem \ref{fcn6}).  
\end{rem}

\begin{cor}
	If $T$ is a tree of order $n\geq 4$ of the form corona of a tree with $K_1$, i.e., $T_1\circ K_1$, where $T_1$ is a tree, then  $$\mathcal{C}_f(T)\leq \frac{n}{2}+2.$$
\end{cor}
\begin{proof}
	It suffices to show that $\gamma_f(T_1\circ K_1)=\frac{n}{2}$. We know that 
	$\gamma(T)=\frac{n}{2}$ and since for any graph $G$, $\gamma(G)\leq \gamma_f(G)$, so $\gamma_f(T)\geq \frac{n}{2}$. On the other hand the set of leaves of $T_1\circ K_1$ form a 1-FD set and so $\gamma_f(T)\leq \frac{n}{2}$. Therefore we have the result by  Theorem \ref{upper}.
\end{proof}

\subsection{Fair coalition of paths and cycles}

This subsection focuses on determining the fair coalition number for paths and cycles, starting with paths. We need the following Lemma:

\begin{lem}{\rm\cite{6}}\label{cnpath}
	The coalition number of any path $P_n$ is at most $6$, i.e., $C(P_n)\leq 6$.
\end{lem}

First we state the following 
observation for path $P_n$ with $V(P_n)=\{v_1,v_2,\dots,v_n\}$ and $E(P_n)=\{\{v_1,v_2\},\{v_2,v_3\},\dots,\{v_{n-1},v_n\}\}$.

\begin{obs}
	\begin{enumerate}
		\item[(i)] $C_f(P_{4})=C_f(P_{5})=C_f(P_6)=C_f(P_8)=4$.
		\item[(ii)] $C_f(P_{7})=C_f(P_9)=C_f(P_{10})=C_f(P_{11})=C_f(P_{13})=5$.
		\item[(iii)] $C_f(P_{ 12})=6.$
	\end{enumerate} 
\end{obs}
\proof 
Note that for any graph $C_f(G)\leq C(G)$ and so by Lemma \ref{cnpath}, the fair coalition number of any path $P_n$ is at most $6$.
In the following parts we present $fc$-partition with maximum size. 
\begin{enumerate}
	\item[(i)]
	For path $P_4$ the  members of the $fc$-partition with maximum size are $A_1=\lbrace v_1 \rbrace$, $A_2=\lbrace v_2\rbrace$, $A_3=\lbrace v_3\rbrace$ and $A_4=\lbrace v_4\rbrace$. Note that $A_1, A_4$ and $ A_2, A_3$ are partners. 
	
	It is easy to observe that there is no $fc$-partition for $P_5$ with cardinality $5$. For $P_5$  the  members of the $fc$-partition with maximum size are $A_1=\lbrace v_1, v_5 \rbrace$, $A_2=\lbrace v_2\rbrace$, $A_3=\lbrace v_3\rbrace$ and $A_4=\lbrace v_4\rbrace$. Note that $A_1, A_2$ and $ A_1, A_3$ are partners. Also $A_1, A_4$ are partners.
	
	We can see, there is no $fc$-partition for $P_6$ with cardinality $5$ or $6$.
	For $P_6$  the  members of  the $fc$-partition with maximum size are $A_1=\lbrace v_1, v_5 \rbrace$, $A_2=\lbrace v_2, v_6\rbrace$, $A_3=\lbrace v_3\rbrace$ and $A_4=\lbrace v_4\rbrace$. Note that $A_1, A_2$ and $ A_1, A_4$ are partners. Also $A_2, A_3$ are partners.

	First we show that there is no $fc$-partition of size $5$ or $6$ for $P_8$. First suppose on contrary that $C_f(P_8)=6$. So any $fc$-partition consists of two sets of size two and four singleton sets, where each singleton set forms a fair coalition with a set of size two, and each set of size two forms a fair coalition with exactly two singleton sets. Since the vertices $v_3$ and $v_6$ are not in any such fair coalitions, so we have a contradiction. Now we show that $C_f(P_8)\neq 5$. Suppose on contrary that $C_f(P_8)=5$ and $\pi$ is a $fc$-partition of size $5$. We have two cases for the $\pi$: 
	
	\noindent Case 1: The partition $\pi$ contains a set of size three, a set of size two and three singleton sets. 
	
	We know that 
	$\gamma_f(P_8)=3$ and so each singleton set cannot form a $fc$-coalition with a non-singleton set. 
	There exists two singleton sets $A_1$ and $A_2$ in $fc$-partition $\pi$ such that $A_1$ forms a $fc$-coalition with the set of size two  and $A_2$ forms a $fc$-coalition with the set of size three. But this impossible, because the graph $P_8$ has six $1$-fair dominating sets of order $4$ two $1$-fair dominating sets of order three. So we have the contradiction. 
	
	\noindent Case 2: The partition $\pi$ contains three sets of size two and two singleton sets.
	
	Since the graph $P_8$ has exactly two fair dominating sets of order three, so this case does not happen.
	
	For the path $P_8$  the  members of the $fc$-partition with maximum size are $A_1=\lbrace v_1, v_8 \rbrace$, $A_2=\lbrace v_2, v_7\rbrace$, $A_3=\lbrace v_3, v_6 \rbrace$ and $A_4=\lbrace v_4, v_5\rbrace$. Note that $ A_1, A_4$ and $A_2, A_3$ are partners.

	\item[(ii)] 
	It is easy to see that, there is no $fc$-partition for $P_6$ with cardinality $6$.
	The partition $\{\{v_2,v_6\},\{v_1,v_7\},\{v_3\},\{v_4\},\{v_5\}\}$ is a $fc$-partition of maximum size for $P_7$, and so $C_f(P_7)=5$.
	
	Now we consider $P_9$. 
	We show that $C_f(P_9)\neq 6$. Suppose on contrary that $C_f(P_9)= 6$. We consider the following cases: 
	
	\noindent Case 1: 
	The $fc$-partition $\pi$ contains a set of size three, a set of size two and four singleton sets. Since $\gamma_f(P_9)=3$, each singleton set must form a fair  coalition with a non-singleton set and so each non-singleton set must form a fair coalition with two singleton sets. 
	But we can see that the graph $P_9$ only has one a fair dominating set of size three, which is $\{v_2,v_5,v_8\}$. Therefore we have a contradiction.
	
	\noindent Case 2:  The $fc$-partition $\pi$ contains  three sets of size two and three singleton sets. Similar to Case 1, we can obtain a contradiction. 
	
	The partition $\{\{v_3,v_6,v_9\},\{v_1,v_4,v_8\},\{v_2\},\{v_5\},\{v_7\}\}$ is a fair partition for $P_9$, so $C_f(P_9) =5$.
	
	The proof for $P_{10},P_{11}$ and $P_{13}$ are similar. 
	
	\item[(iii)] The proof is similar to the proof of Theorem \ref{fcpath6}.
	
\end{enumerate}

\begin{thm}\label{fcpath6}
	For $n \geq 14$, $C_f(P_{n})=6$.
\end{thm}

\proof
Suppose that $V(P_n)=\lbrace v_1, v_2,\dots,v_n \rbrace.$ By Lemma \ref{cnpath} it is enough to present $fc$-partition with cardinality six. Let consider two cases:\\
Case 1) If $n$ is even, then for $k\geq 3$ the $fc$-partition with maximum size is as follows:
$$\Upsilon= \lbrace A_1, A_2, A_3, A_4, A_5, A_6\rbrace,$$
where $A_1=\lbrace v_1, v_4, v_7, v_{4k-1}, \dots, v_{4k+2}, \dots \rbrace$, $A_2=\lbrace v_2, v_6, v_9, v_{4k},\dots, v_{4k+1}, \dots \rbrace$, $A_3=\lbrace v_3\rbrace$, $A_4=\lbrace v_{10}\rbrace$, $A_5=\lbrace v_5\rbrace$ and $A_6=\lbrace v_8\rbrace$.\\
Note that $A_1, A_4$ and $A_1, A_6$ are partners. Also $A_2, A_3$ and $A_2, A_5$ are partners.\\
Case 2) If $n$ is odd, then for $k\geq 4$ the $fc$-partition with maximum size is as follows:
$$\Upsilon= \lbrace A_1, A_2, A_3, A_4, A_5, A_6\rbrace$$
Where $A_1=\lbrace v_1, v_4, v_7, v_{10}, v_{4k-2}, \dots, v_{4k+1}, \dots \rbrace$, $A_2=\lbrace v_2,v_9, v_{12}, v_{4k-1}, \dots, v_{4k}, \dots \rbrace$, $A_3=\lbrace v_3, v_6\rbrace$ and $A_4=\lbrace v_{13}\rbrace$, $A_5=\lbrace v_5, v_8\rbrace$ and $A_6=\lbrace v_{11}\rbrace$.\\
Note that $A_1, A_4$ and $A_1, A_6$ are partners. Also $A_2, A_3$ and $A_2, A_5$ are partners. \qed

Here, we compute the fair coalition number of cycle $C_n$.
We need the following lemme:

\begin{lem}\label{cncycle}
	The fair  coalition number of any cycle $C_n$ is at most $6$, i.e., $C_f(C_n)\leq 6$.
\end{lem}
\proof 
We know that (\cite{6}) the fair  coalition number of any cycle $C_n$ is at most $6$, and since $C_f(C_n)\leq C(C_n)$, so we have the result.

\begin{thm}
	For $n\geq 10$, $(n\neq 11)$, $C_f(C_n)=6$.
\end{thm}
\begin{proof}
	Suppose that $V(C_n)=\{v_1,v_2,...,v_n\}$. By Lemma \ref{cncycle}, it is sufficient to present a $fc$-partition of size six. 
	\begin{enumerate}  
		\item[(i)] 
		If $n=3k$, we present $fc$-partition with maximum size. Let consider two cases:
		
		Case 1) If $k$ is odd, then the $fc$-partition with maximum size is as follows:
		$$\Upsilon=\{A_1,A_2,A_3,A_4,A_5, A_6\},$$
		where 
		$A_1=\{v_1,v_4,...,v_{\lfloor\frac{3k}{2}\rfloor}\}$, $A_2=\{v_2,v_5,...,v_{\lfloor\frac{3k}{2}\rfloor+1}\},$
		$A_3=\{v_3,v_6,...,v_{\lfloor\frac{3k}{2}\rfloor+2}\},$ \\
		$A_4=\{v_{\lfloor\frac{3k}{2}\rfloor+3},v_{\lfloor\frac{3k}{2}\rfloor+6},...,v_{3k-2}\}$,
		$A_5=\{v_{\lfloor\frac{3k}{2}\rfloor+4},v_{\lfloor\frac{3k}{2}\rfloor+7},...,v_{3k-1}\}$,\\
		$A_6=\{v_{\lfloor\frac{3k}{2}\rfloor+5},v_{\lfloor\frac{3k}{2}\rfloor+8},...,v_{3k}\}$.
		
		Note that $A_1$, $A_4$  and $A_2$, $A_5$ are partner. Also $A_3$ and $A_6$ are partner.

		Case 2) If $k$ is even, then the $fc$-partition with maximum size is as follows:
		$$\Upsilon=\{A_1,A_2,A_3,A_4,A_5, A_6\},$$
		where 
		$A_1=\{v_1,v_4,...,v_{\lfloor\frac{3k}{2}\rfloor-2}\}$, $A_2=\{v_2,v_5,...,v_{\lfloor\frac{3k}{2}\rfloor-1}\},$
		$A_3=\{v_3,v_6,...,v_{\lfloor\frac{3k}{2}\rfloor}\},$ \\
		$A_4=\{v_{\lfloor\frac{3k}{2}\rfloor+1},v_{\lfloor\frac{3k}{2}\rfloor+4},...,v_{3k-2}\}$,
		$A_5=\{v_{\lfloor\frac{3k}{2}\rfloor+2},v_{\lfloor\frac{3k}{2}\rfloor+5},...,v_{3k-1}\}$,\\
		$A_6=\{v_{\lfloor\frac{3k}{2}\rfloor+3},v_{\lfloor\frac{3k}{2}\rfloor+6},...,v_{3k}\}$.
		
		Note that $A_1$, $A_4$  and $A_2$, $A_5$ are partner. Also $A_3$ and $A_6$ are partner. Therefore $C_f(C_{3k})=6$. 
		
		\item[(ii)] 
		If $n=3k+1$, the $fc$-partition of $C_{3k+1}$ with maximum size is as follows:
		$$\Upsilon=\{A_1,A_2,A_3,A_4,A_5, A_6\},$$
		where for odd $k$ we have;\\
		$A_1=\{v_i \mid  i=\lbrace 1, 4,..., 3k-2\rbrace \}$, $A_2=\{v_{3k-1}\},$
		$A_3=\{v_{3k+1}\},$\\
		$A_4=\{v_2, v_{3k-3}, v_{3k}\}$,
		$A_5=\{v_i\mid  i=\lbrace 3, 6, ..., 3k-6\rbrace\}$,\\
		$A_6=\{ v_i\mid i=\lbrace 5, 8, ..., 3k-4\rbrace \}.$
		
		And for even $k$;\\
		$A_1=\{v_i \mid  i=\lbrace 1, 4,8, 11, 14,..., 3k-1\rbrace \}$, $A_2=\{v_5 \},$
		$A_3=\{v_7 \},$\\
		$A_4=\{v_3, v_6, v_9,..., v_{3k-3}\}$,
		$A_5=\{v_2, v_{3k}\}$,\\
		$A_6=\{ v_{3k-2}, v_{3k+1} \}.$

		Note that $A_1$, $A_2$  and $A_1$, $A_3$ are partner. Also $A_4$, $A_5$  and $A_4$, $A_6$ are partner. So $C_f(C_{3k+1})=6$.
		
		\item[(iii)] 
		The $fc$-partition of $C_{3k+2}$ with maximum size is as follows:
		$$\Upsilon=\{A_1,A_2,A_3,A_4, A_5, A_6\},$$
		where for even $k$;\\
		$A_1=\{v_i \mid  i=\lbrace 1, 4,8, 11, 14,..., 3k+2\rbrace \}$, $A_2=\{v_5\},$
		$A_3=\{v_7\},$ \\
		$A_4=\{v_2, v_3, v_6, v_9, v_{13}, v_{16}, ..., v_{3k+1}\}$, $A_5=\{ v_{10} \}$, $A_6=\{ v_{12} \}.$
		
		And for odd $k$;\\
		$A_1=\lbrace v_1, v_2, v_5, v_8, ..., v_{ 3k-7}, v_{3k} \rbrace$, $A_2=\lbrace v_{3k-4}, v_{3k-1}\rbrace,$\\
		$A_3=\lbrace v_{3k-6}, v_{3k-3} \rbrace$, $A_4=\lbrace v_3, v_7, v_{10}, v_{13}, v_{16}, ..., v_{3k+1}, v_{3k+2}\rbrace$, $A_5=\lbrace v_4 \rbrace $ and $A_6=\lbrace v_6 \rbrace.$
		
		Note that $A_1$, $A_2$  and $A_1$, $A_3$ are partner. Also $A_4$, $A_5$ and $A_4$, $A_6$ are partner. So $C_f(C_{3k+2})=6$.		  
		
	\end{enumerate} 	
\end{proof}

\section{Fair coalition of cubic graphs of order at most $10$}

In this section, we obtain the fair coalition number of cubic graphs of order at most $10$.  In particular, we obtain the fair coalition number of the Petersen graph. 
The coalition number and the total coalition number of cubic graphs of order at most $10$ have studied in \cite{1} and \cite{CCO2}, respectively. 
We need the following lemma:

\begin{lem}{\rm\cite{1}}
	If $G$ is a cubic graph of 
	order at most $10$, then $C(G)\in \{6, 7, 8\}$.  
\end{lem}
So we have the following corollary:

\begin{cor}\label{fcg10}
	If $G$ is a cubic graph of 
	order at most $10$, then $C_f(G)\leq 8$.  
\end{cor}

\subsection{Results for cubic graphs of order $6$}
\begin{figure}[h!]
	\hglue2.5cm
	\includegraphics[width=9cm,height=3.3cm]{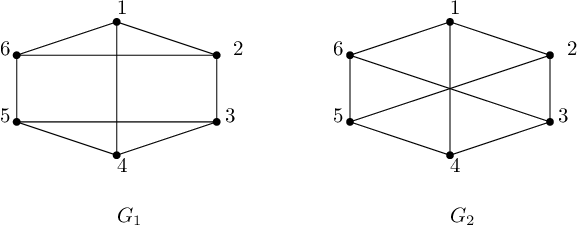}
	\vglue-10pt \caption{\label{Cubic6} Cubic graphs of order $6$.}
\end{figure}

In this subsection, we obtain the fair coalition number of cubic graphs of order $6$.  There are exactly two  cubic graphs of order $6$ which are denoted by $G_{1}$ and $G_{2}$ in Figure \ref{Cubic6}.

\begin{thm}\label{fcn6}
	The fair coalition number of cubic graphs $G_1$ and $G_2$ of order $6$ is $6$.
\end{thm}

\proof
Suppose that $V(G_1)=\lbrace 1, 2, 3, 4, 5, 6 \rbrace$ and $ V(G_2)=\lbrace 1, 2, 3, 4, 5, 6 \rbrace.$\\
By Corollary \ref{fcg10}, it is sufficient to present the $fc$-partition with maximum size as follows;
$$\Upsilon= \lbrace A_1, A_2, A_3, A_4, A_5, A_6\rbrace.$$
Where $A_1=\lbrace 1 \rbrace$, $A_2=\lbrace 2 \rbrace$, $A_3=\lbrace 3 \rbrace$, $A_4=\lbrace 4 \rbrace$, $A_5=\lbrace 5 \rbrace$ and $A_6=\lbrace 6 \rbrace$. 
Note that $A_1, A_4$ and $A_2, A_3$ are partners. Also $A_5, A_6$ are partners. \qed

\subsection{Results for cubic graphs of order $8$}

\begin{figure}[h!]
	\hglue1.5cm 
	\includegraphics[width=11cm,height=7.2cm]{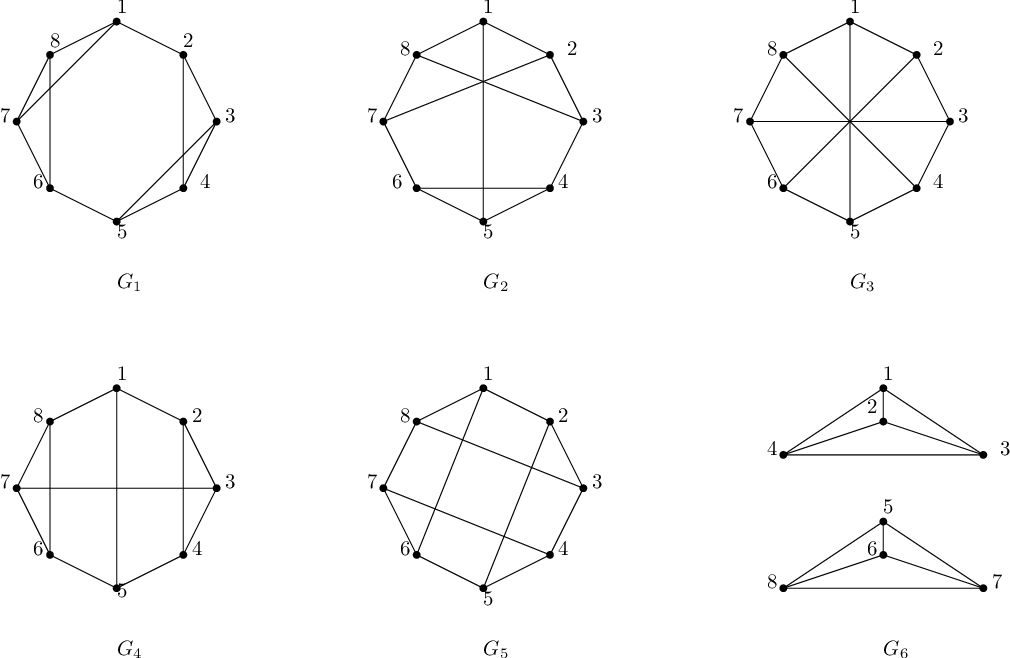}
	\vglue-10pt \caption{\label{Cubic8} Cubic graphs of order $8$.}
\end{figure}

In the following we obtain the fair coalition number of cubic graphs of order $8$. 
There are exactly $6$ cubic graphs of order $8$ which are denoted by  $G_{1},G_{2},...,G_{6}$ in Figure \ref{Cubic8}. 

\begin{thm}
	\begin{enumerate} 
		\item[(i)] 
		For the cubic graph $G_1$ of order $8$, $C_f(G_1)=8$.
		\item[(ii)] 
		For the cubic graph $G_2$ of order $8$, $C_f(G_2)=5$.
		\item[(iii)] 
		For the cubic graph $G_3$ of order $8$, $C_f(G_3)=5$.
		\item[(iv)] 
		For the cubic graph $G_4$ of order $8$, $C_f(G_4)=6$.	
		\item[(v)] 
		For the cubic graph $G_5$ of order $8$, $C_f(G_5)=8$.
		\item[(vi)] 
		For the cubic graph $G_6$ of order $8$, $C_f(G_6)=8$.
	\end{enumerate} 	
\end{thm}

\proof
Consider the cubic graphs of order $8$ in Figure \ref{Cubic8}. 
\begin{enumerate}
	
	\item[(i)] 
	By Corollary \ref{fcg10}, $C_f(G_1)\leq 8$. So it is enough to present the $fc$-partition with size $8$ as follows;
	$$\Upsilon= \lbrace A_1, A_2, A_3, A_4, A_5, A_6, A_7, A_8\rbrace,$$
	where $A_1=\lbrace 1 \rbrace$, $A_2=\lbrace 2 \rbrace$, $A_3=\lbrace 3 \rbrace$, $A_4=\lbrace 4 \rbrace$, $A_5=\lbrace 5 \rbrace$, $A_6=\lbrace 6 \rbrace$, $A_7=\lbrace 7 \rbrace$ and $A_8=\lbrace 8 \rbrace$. 
	Note that $A_1, A_5$ and $A_2, A_6$ are partners. Also $A_3, A_7$ and $A_4, A_8$ are partners. 
	
	\item[(ii)] 
	It is easy to see that there is no $fc$-partition of $G_2$ of size  $6,7,8$. 
	We present $fc$-partition with size $5$ as follows;
	$$\Upsilon= \lbrace A_1, A_2, A_3, A_4, A_5\rbrace,$$
	where $A_1=\lbrace 1, 5\rbrace$, $A_2=\lbrace 2 \rbrace$, $A_3=\lbrace 3, 4 \rbrace$, $A_4=\lbrace 6, 7 \rbrace$ and 
	$A_5=\lbrace 8 \rbrace$. 
	Note that $A_1, A_2$ and $A_2, A_3$ are partners. Also $A_4, A_5$ are partners. 
	
	\item[(iii)] 
	It is easy to see that there is no $fc$-partition of $G_3$ of size  $6,7,8$. 
	We present $fc$-partition with size $5$ as follows;
	$$\Upsilon= \lbrace A_1, A_2, A_3, A_4, A_5\rbrace,$$
	where $A_1=\lbrace 1 \rbrace$, $A_2=\lbrace 2, 3 \rbrace$, $A_3=\lbrace 7, 8 \rbrace$, $A_4=\lbrace 4 \rbrace$ 
	and $A_5=\lbrace 5, 6 \rbrace$.
	Note that $A_1, A_2$ and $A_1, A_3$ are partners. Also $A_4, A_5$ are partners. 
	
	\item[(iv)] 
	It is easy to see that there is no $fc$-partition of $G_4$ of size  $7,8$. 
	We present $fc$-partition with size $6$ as follows;
	$$\Upsilon= \lbrace A_1, A_2, A_3, A_4, A_5, A_6\rbrace,$$
	where $A_1=\lbrace 2 \rbrace$, $A_2=\lbrace 6 \rbrace$, $A_3=\lbrace 4 \rbrace$, $A_4=\lbrace 8 \rbrace$, $A_5=\lbrace 1, 5 \rbrace$ and $A_6=\lbrace 3, 7 \rbrace$.
	Note that $A_1, A_2$ and $A_3, A_4$ are partners. Also $A_5, A_6$ are partners. 
	
	\item[(v)] 
	We present $fc$-partition with maximum size as follows;
	$$\Upsilon= \lbrace A_1, A_2, A_3, A_4, A_5, A_6, A_7, A_8\rbrace,$$
	where $A_1=\lbrace 1 \rbrace$, $A_2=\lbrace 2 \rbrace$, $A_3=\lbrace 3 \rbrace$, $A_4=\lbrace 4 \rbrace$, $A_5=\lbrace 5 \rbrace$, $A_6=\lbrace 6 \rbrace$, $A_7=\lbrace 7 \rbrace$ and $A_8=\lbrace 8 \rbrace$. 
	Note that $A_1, A_4$ and $A_2, A_7$ are partners. Also $A_3, A_6$ and $A_5, A_8$ are partners. 
	
	\item[(vi)] 
	We present $fc$-partition with maximum size as follows;
	$$\Upsilon= \lbrace A_1, A_2, A_3, A_4, A_5, A_6, A_7, A_8\rbrace,$$
	where $A_1=\lbrace 1 \rbrace$, $A_2=\lbrace 2 \rbrace$, $A_3=\lbrace 3 \rbrace$, $A_4=\lbrace 4 \rbrace$, $A_5=\lbrace 5 \rbrace$, $A_6=\lbrace 6 \rbrace$, $A_7=\lbrace 7 \rbrace$ and $A_8=\lbrace 8 \rbrace$. 
	Note that $A_1, A_5$ and $A_2, A_6$ are partners. Also $A_3, A_7$ and $A_4, A_8$ are partners. 
	\qed
	
\end{enumerate}

\subsection{Results for cubic graphs of order $10$} 

\begin{figure}[h!]
	\hglue1.25cm
	\includegraphics[width=11cm,height=4.89cm]{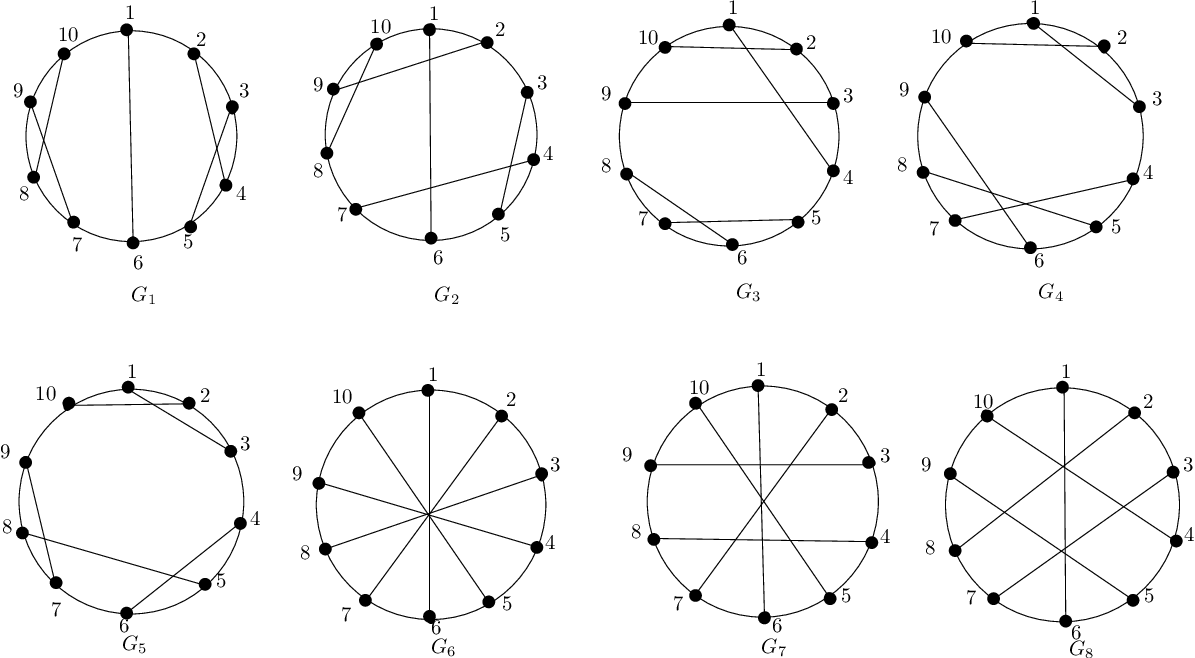}
	\vglue5pt
	\hglue1.25cm
	\includegraphics[width=11cm,height=4.99cm]{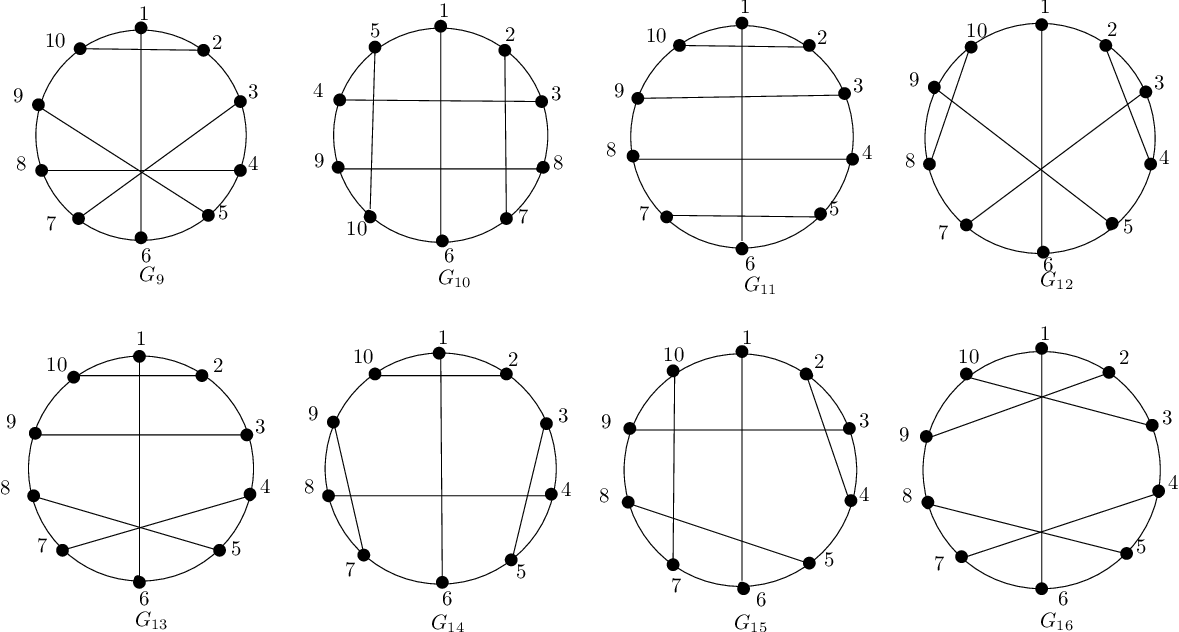}
	\hglue1.25cm
	\includegraphics[width=10.7cm,height=4.9cm]{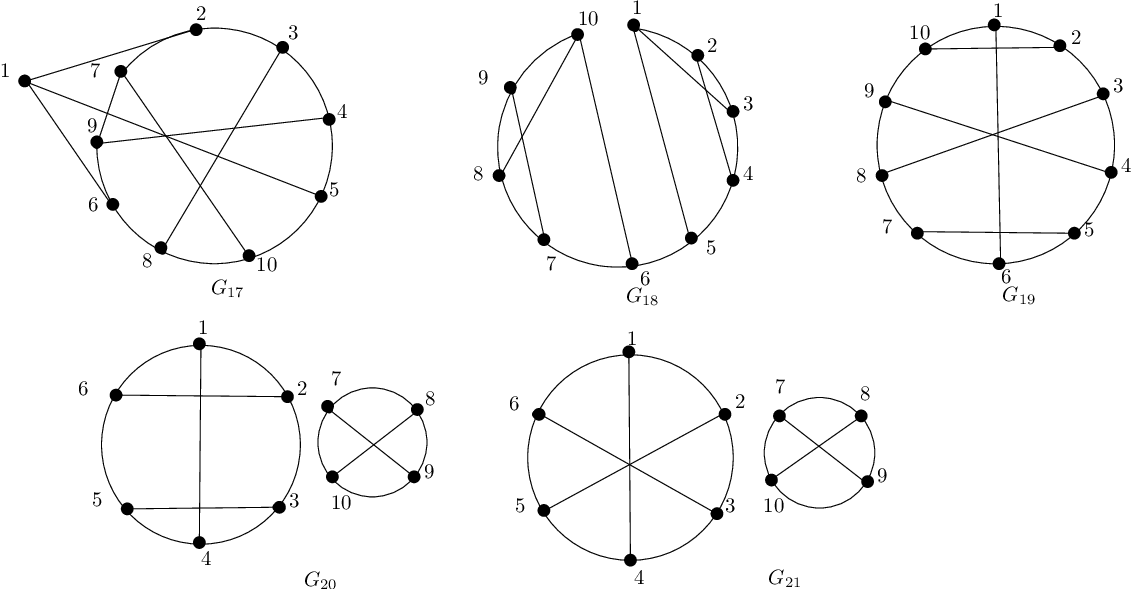}
	\hglue1.25cm
	\vglue-10pt \caption{\label{figure2} Cubic graphs of order $10$.}
\end{figure}

In this subsection, we obtain the fair coalition number of cubic graphs of order $10$.  There are exactly $21$ cubic graphs of order $10$ denoted by  $G_{1},G_{2},...,G_{21}$ in Figure \ref{figure2} (see \cite{1}). 
In particular, the graph $G_{17}$ is isomorphic to the Petersen graph $P$.

Now we state and prove the following theorem.

\begin{thm}
	Let $G_i$ ($1\leq i\leq 21$) be the cubic graphs of order $10$. Then 
	$C_f(G_i)=4$ for $i \in \lbrace 1, 12, 14, 17, 18, 19\rbrace.$
\end{thm}

\proof Consider the cubic graphs $G_1, G_2, \cdots, G_{21}$ of order $10$ as shown in Figure \ref{figure2}. 

It is easy to see that, there is no $fc$-partition for $G_1$ of size $5,6,7,8$.
We present $fc$-partition with maximum size as follows for $G_1$:
$$\Upsilon= \lbrace A_1, A_2, A_3, A_4\rbrace,$$
where $A_1=\lbrace 1,6\rbrace$, $A_2=\lbrace 2,{10}\rbrace$, $A_3=\lbrace 5, 7\rbrace$ and $A_4=\lbrace 3, 4, 8, 9\rbrace$. 
Note that $A_1, A_2$ and $A_1, A_3$ are partners. Also $A_1, A_4$ are partners. 

\medskip
There is no $fc$-partition for $G_{12}$ of size $5,6,7,8$.
The following partition is the $fc$-partition with maximum size for $G_{12}$. 
$$\Upsilon= \lbrace A_1, A_2, A_3, A_4\rbrace,$$
where $A_1=\lbrace 1, 6\rbrace$, $A_2=\lbrace 2, 10\rbrace$, $A_3=\lbrace 3,4, 8, 9\rbrace$ and $A_4=\lbrace 5, 7\rbrace$.
Note that $A_1, A_2$ and $A_1, A_3$ are partners. Also $A_1, A_4$ are partners.

\medskip 
There is no $fc$-partition for $G_{14}$ of size $5,6,7,8$.
The following partition is the $fc$-partition with maximum size for $G_{14}$. 
$$\Upsilon= \lbrace A_1, A_2, A_3, A_4\rbrace,$$
where $A_1=\lbrace 1, 4, 5, 9\rbrace$, $A_2=\lbrace 2, 7\rbrace$, $A_3=\lbrace 3,6\rbrace$ and $A_4=\lbrace 8, {10}\rbrace$. 
Note that $A_1, A_2$ and $A_1, A_4$ are partners. Also $A_3, A_4$ are partners.

\medskip 
There is no $fc$-partition for $G_{17}$ of size $5,6,7,8$.
The following partition is the $fc$-partition with maximum size for $G_{17}$. 
$$\Upsilon= \lbrace A_1, A_2, A_3, A_4\rbrace,$$
where $A_1=\lbrace 1, 3, 8\rbrace$, $A_2=\lbrace 2, 6\rbrace$, $A_3=\lbrace 4,5, {10}\rbrace$ and $A_4=\lbrace 7, 9\rbrace$.  
Note that $A_1, A_2$ and $A_3, A_4$ are partners.

\medskip
There is no $fc$-partition for $G_{18}$ of size $5,6,7,8$.
The following partition is the $fc$-partition with maximum size for $G_{18}$.
$$\Upsilon= \lbrace A_1, A_2, A_3, A_4\rbrace,$$
where $A_1=\lbrace 1, 7\rbrace$, $A_2=\lbrace 4, {10}\rbrace$, $A_3=\lbrace 5, 6\rbrace$ and $A_4=\lbrace 2, 3, 8, 9\rbrace$. 
Note that $A_1, A_2$ and $A_1, A_3$ are partners. Also $A_3, A_4$ are partners.

\medskip 
There is no $fc$-partition for $G_{19}$ of size $5,6,7,8$.
The following partition is the $fc$-partition with maximum size for $G_{19}$.
$$\Upsilon= \lbrace A_1, A_2, A_3, A_4\rbrace,$$
where $A_1=\lbrace 1, 8, 9\rbrace$, $A_2=\lbrace 2\rbrace$, $A_3=\lbrace 3, 4, 5 \rbrace$ and $A_4=\lbrace 6, 7, {10}\rbrace$. 
Note that $A_1, A_3$ and $A_2, A_3$ are partners. Also $A_3, A_4$ are partners. \qed

\begin{thm}
	Let $G_i$ ($1\leq i\leq 21$) be the cubic graphs of order $10$. Then 
	$C_f(G_i)=5$ for $i \in \lbrace 2, 6, 7, 8, 9, 11, 13, 16\rbrace.$
\end{thm}
\proof
It is easy to see that there is no $fc$-partition for $G_{2}$ of size $6,7,8$.
The following partition is the $fc$-partition with maximum size for $G_{2}$.
$$\Upsilon= \lbrace A_1, A_2, A_3, A_4, A_5\rbrace,$$
where $A_1=\lbrace 1, 6\rbrace$, $A_2=\lbrace 2, 7\rbrace$, $A_3=\lbrace 3, 4\rbrace$, $A_4=\lbrace 5, {10} \rbrace$, $A_5=\lbrace  8, 9 \rbrace$. 
Note that $A_1, A_2$ and $A_3, A_4$ are partners. Also $A_4, A_5$ are partners. 

\medskip 
There is no $fc$-partition for $G_{6}$ of size $6,7,8$.
The following partition is the $fc$-partition with maximum size for $G_{6}$.
$$\Upsilon= \lbrace A_1, A_2, A_3, A_4, A_5\rbrace,$$
where $A_1=\lbrace 1, 2\rbrace$, $A_2=\lbrace 3, 4 \rbrace$, $A_3=\lbrace 5, 6\rbrace$, $A_4=\lbrace 7, 8 \rbrace$, $A_5=\lbrace 9, {10}\rbrace$.  
Note that $A_1, A_2$ and $A_2, A_3$ are partners. Also $A_4, A_5$ are partners.

\medskip 
There is no $fc$-partition for $G_{7}$ of size $6,7,8$.
The following partition is the $fc$-partition with maximum size for $G_{7}$.
$$\Upsilon= \lbrace A_1, A_2, A_3, A_4, A_5\rbrace,$$
where $A_1=\lbrace 1, 2\rbrace$, $A_2=\lbrace 3, 4\rbrace$, $A_3=\lbrace 5, 6\rbrace$, $A_4=\lbrace 7, 8\rbrace$, $A_5=\lbrace 9, {10}\rbrace$.  
Note that $A_1, A_2$ and $A_2, A_3$ are partners. Also $A_4, A_5$ are partners.

\medskip 
There is no $fc$-partition for $G_{8}$ of size $6,7,8$.
The following partition is the $fc$-partition with maximum size for $G_{8}$.
$$\Upsilon= \lbrace A_1, A_2, A_3, A_4, A_5\rbrace,$$
where $A_1=\lbrace 1, 2\rbrace$, $A_2=\lbrace 5, 6\rbrace$, $A_3=\lbrace 3, 4\rbrace$, $A_4=\lbrace 7, {10}\rbrace$, $A_5=\lbrace 8, 9\rbrace$.
Note that $A_1, A_2$ and $A_3, A_4$ are partners. Also $A_4, A_5$ are partners.

\medskip 
There is no $fc$-partition for $G_{9}$ of size $6,7,8$.
The following partition is the $fc$-partition with maximum size for $G_{9}$.
$$\Upsilon= \lbrace A_1, A_2, A_3, A_4, A_5\rbrace,$$
where $A_1=\lbrace 1, 6\rbrace$, $A_2=\lbrace 2, 3 \rbrace$, $A_3=\lbrace 4, 5\rbrace$, $A_4=\lbrace 7, 8 \rbrace$, $A_5=\lbrace 9, {10}\rbrace$.  
Note that $A_1, A_3$ or $A_1, A_4$ are partners. Also $A_2, A_3$ and $A_4, A_5$ are partners.

\medskip 
There is no $fc$-partition for $G_{11}$ of size $6,7,8$.
The following partition is the $fc$-partition with maximum size for $G_{11}$.
$$\Upsilon= \lbrace A_1, A_2, A_3, A_4, A_5\rbrace,$$
where $A_1=\lbrace 1\rbrace$, $A_2=\lbrace 2, 5, 7, {10}\rbrace$, $A_3=\lbrace 3, 9\rbrace$, $A_4=\lbrace 4, 8 \rbrace$, $A_5=\lbrace 6\rbrace$. 
Note that $A_1, A_4$ and $A_2, A_4$ or $A_2, A_3$ are partners. Also $A_3, A_5$ are partners.

\medskip 
There is no $fc$-partition for $G_{13}$ of size $6,7,8$.
The following partition is the $fc$-partition with maximum size for $G_{13}$.
$$\Upsilon= \lbrace A_1, A_2, A_3, A_4, A_5\rbrace,$$
where $A_1=\lbrace 1, 2, 5, 8\rbrace$, $A_2=\lbrace 3, 6\rbrace$, $A_3=\lbrace 4\rbrace$, $A_4=\lbrace 7, {10} \rbrace$, $A_5=\lbrace 9\rbrace$. 
Note that $A_1, A_2$ and $A_2, A_5$ are partners. Also $A_3, A_4$ are partners.

\medskip 
There is no $fc$-partition for $G_{16}$ of size $6,7,8$.
The following partition is the $fc$-partition with maximum size for $G_{16}$.
$$\Upsilon= \lbrace A_1, A_2, A_3, A_4, A_5\rbrace,$$
where $A_1=\lbrace 1, 2 \rbrace$, $A_2=\lbrace 6, 7\rbrace$, $A_3=\lbrace 3, 4\rbrace$, $A_4=\lbrace 5, {10} \rbrace$, $A_5=\lbrace 8, 9\rbrace$.  
Note that $A_1, A_2$ and $A_3, A_4$ are partners. Also $A_4, A_5$ are partners.\qed

\begin{thm}
	Let $G_i$ ($1\leq i\leq 21$) be the cubic graphs of order $10$. Then 
	$C_f(G_i)=7$ for $i \in \lbrace 3, 4, 5, 10, 15, 20, 21\rbrace.$
\end{thm}
\proof 
It is easy to see that there is no $fc$-partition of size $8$ for $G_i$, where $i \in \lbrace 3, 4, 5, 10, 15, 20, 21\rbrace$. So for $i \in \lbrace 3, 4, 5, 10, 15, 20, 21\rbrace$
we present the $fc$-partition of  size $7$.

The following partition is the $fc$-partition with maximum size for $G_{3}$.
$$\Upsilon= \lbrace A_1, A_2, A_3, A_4, A_5, A_6, A_7\rbrace,$$
where $A_1=\lbrace 1, 4\rbrace$, $A_2=\lbrace 8\rbrace$, $A_3=\lbrace 2, 3\rbrace$, $A_4=\lbrace 6 \rbrace$, $A_5=\lbrace 5\rbrace$, $A_6=\lbrace 9, {10}\rbrace$, $A_7=\lbrace 7 \rbrace$. 
Note that $A_1, A_2$ and $A_3, A_4$ are partners. Also $A_5, A_6$  and $A_3, A_7$ are partners. 

\medskip 
The following partition is the $fc$-partition with maximum size for $G_{4}$.
$$\Upsilon= \lbrace A_1, A_2, A_3, A_4, A_5, A_6, A_7\rbrace,$$
where $A_1=\lbrace 1, 5\rbrace$, $A_2=\lbrace 2, 6\rbrace$, $A_3=\lbrace 3\rbrace$, $A_4=\lbrace 8\rbrace$, $A_5= \lbrace 7 \rbrace$, $A_6=\lbrace 4, 9 \rbrace$ and $A_7=\lbrace {10}\rbrace$.  
Note that $A_1, A_4$ and $A_2, A_5$ are partners. Also $A_3, A_6$ and $A_6, A_7$ are partners.

\medskip 
The following partition is the $fc$-partition with maximum size for $G_{5}$.
$$\Upsilon= \lbrace A_1, A_2, A_3, A_4, A_5, A_6, A_7\rbrace,$$
where $A_1=\lbrace 1, 5\rbrace$, $A_2=\lbrace 2, 6\rbrace$, $A_3=\lbrace 3\rbrace$, $A_4=\lbrace 8\rbrace$, $A_5= \lbrace 7 \rbrace$, $A_6=\lbrace 4, 9 \rbrace$ and $A_7=\lbrace {10}\rbrace$.  
Note that $A_1, A_4$ and $A_2, A_5$ are partners. Also $A_3, A_6$ and $A_6, A_7$ are partners.

\medskip 
The following partition is the $fc$-partition with maximum size for $G_{10}$.
$$\Upsilon= \lbrace A_1, A_2, A_3, A_4, A_5, A_6, A_7\rbrace,$$
where $A_1=\lbrace 1\rbrace$, $A_2=\lbrace 2\rbrace$, $A_3=\lbrace 3, 6\rbrace$, $A_4=\lbrace 4, 5\rbrace$, $A_5= \lbrace 7, 8 \rbrace$, $A_6=\lbrace  9 \rbrace$ and $A_7=\lbrace {10}\rbrace$. 
Note that $A_1, A_5$ and $A_2, A_5$ are partners. Also $A_3, A_6$ and $A_4, A_7$ are partners.

\medskip 
The following partition is the $fc$-partition with maximum size for $G_{15}$.
$$\Upsilon= \lbrace A_1, A_2, A_3, A_4, A_5, A_6, A_7\rbrace,$$
where $A_1=\lbrace 1,2\rbrace$, $A_2=\lbrace 3, 6\rbrace$, $A_3=\lbrace 4, 5\rbrace$, $A_4=\lbrace  7\rbrace$, $A_5= \lbrace 8 \rbrace$, $A_6=\lbrace  9 \rbrace$ and $A_7=\lbrace {10}\rbrace$.  
Note that $A_1, A_5$ and $A_2, A_4$ are partners. Also $A_2, A_6$ and $A_3, A_7$ are partners.

\medskip 
The following partition is the $fc$-partition with maximum size for $G_{20}$.
$$\Upsilon= \lbrace A_1, A_2, A_3, A_4, A_5, A_6, A_7\rbrace,$$
where $A_1=\lbrace 1,7\rbrace$, $A_2=\lbrace 2, 8\rbrace$, $A_3=\lbrace 3\rbrace$, $A_4=\lbrace  4\rbrace$, $A_5= \lbrace 5, 6 \rbrace$, $A_6=\lbrace  9 \rbrace$ and $A_7=\lbrace {10}\rbrace$. 
Note that $A_1, A_4$ and $A_2, A_3$ are partners. Also $A_5, A_6$ and $A_5, A_7$ are partners.

\medskip 
The following partition is the $fc$-partition with maximum size for $G_{21}$.
$$\Upsilon= \lbrace A_1, A_2, A_3, A_4, A_5, A_6, A_7\rbrace,$$
where $A_1=\lbrace 1,7\rbrace$, $A_2=\lbrace 2, 8\rbrace$, $A_3=\lbrace 3, 6\rbrace$, $A_4=\lbrace  4\rbrace$, $A_5= \lbrace 5 \rbrace$, $A_6=\lbrace  9 \rbrace$ and $A_7=\lbrace {10}\rbrace$.  
Note that $A_1, A_4$ and $A_2, A_5$ are partners. Also $A_3, A_6$ and $A_3, A_7$ are partners. \qed

\section{Conclusion}
This paper introduces the concept of the fair coalition in graphs and explores various properties related to its number. We have demonstrated that when a graph $G$ has at least three vertices without full vertices, then $C_f(G)\geq 2d_f(G)$.  We have determined the precise values of ${C}_f(P_n)$, $\mathcal{C}_f(C_n)$, and the fair coalition number of the cubic graphs of order at most $10$. There is much work to be done in this area.
\begin{enumerate}
	\item What is the fair coalition number of  graph operations, such as corona, Cartesian product, join, lexicographic, and so on?
	
	\item What is the fair coalition number of  natural and fractional powers of a graph?

	\item What is the effects on $\mathcal{C}_f(G)$ when $G$ is modified by operations on vertex and edge of $G$?

	\item Study Nordhaus and Gaddum lower and upper bounds on the sum and the product
	of the fair coalition number of a graph and its complement.
	
	\item Study the complexity of the fair coalition number for many of the graphs.  
	
\end{enumerate}


\begin{thebibliography}{99}
	
	
	
	\bibitem{1} S. Alikhani, H.R. Golmohammadi, E.V. Konstantinova, Coalition of
	cubic graphs of order at most $10$, {\it Comm. Combin.
		Optim.}, {\bf 9(3)}  (2024) 437-450.
	
	
	\bibitem{total}  S. Alikhani, D. Bakhshesh, H. Golmohammadi, Total coalitions in graphs, {\it Quaest. Math.} {\bf 47(11)} (2024)  2283-2294. 
	
	\bibitem{2} S. Alikhani, D. Bakhshesh, H. Golmohammadi, E.V. Konstantinova, Connected
	coalitions in graphs, {\it  Discuss. Math. Graph Theory} {\bf 44} (2024) 1551-1566. 
	
	
	
	\bibitem{Fair1} S. Alikhani, M. Safazadeh, Fair dominating sets of paths, 
	{\it J.  Inform. Optim.  Sci.},   {\bf	44 (5)} (2023) 855-864. 
	
	\bibitem{Fair2} 
	S. Alikhani, M. Safazadeh, On the number of fair dominating sets of graphs, {\it Advanced Studies: Euro
		Tbilisi Math. J.} {\bf  16(1)} (2023) 59-69. 
	
	
	\bibitem{3} D. Bakhshesh, M.A. Henning, D. Pradhan, On the coalition number of trees,  {\it Bull. Malays. Math. Sci. Soc.} (2023) 46:95. 
	
	
	
	
	\bibitem{Henning} Y. Caro, A. Hansberg, M.  Henning,   Fair domination in graphs, {\it Discrete Appl. Math.}  {\bf 312}  (2012) 2905-2914. 
	
	
	\bibitem{5} E.J. Cockayne, S. T. Hedetniemi, Towards a theory of domination in graphs. Net-
	works, 7:247-261, 1977.
	
	\bibitem{CCO2} H.R. Golmohammadi, Total coalitions of cubic graphs of order at most $10$,  {\it Comm. Combin.
		Optim.}, {\bf 10(3)} (2025) 601-615.
	
	
	
	\bibitem{6} T.W. Haynes, J.T. Hedetniemi, S.T. Hedetniemi, A.A. McRae and R. Mohan,
	Introduction to coalitions in graphs, {\it AKCE Int. J. Graphs Combin.}
	{\bf  17(2)} (2020),
	653--659.
	\bibitem{7} T. W. Haynes, J. T. Hedetniemi, S. T. Hedetniemi, A. A. McRae, and R. Mohan, {\it Coalition graphs}, {Commun. Comb. Optim.} {\bf 8(2)} (2023) 423--430.
	
	\bibitem{8} T.W. Haynes, J.T. Hedetniemi, S.T. Hedetniemi, A.A. McRae, R. Mohan,
	Coalition graphs of paths, cycles, and trees, {\it Discuss. Math. Graph Theory}
	{\bf 43} (2023) 931-946.
	
	
	\bibitem{9} T.W. Haynes, J.T. Hedetniemi, S.T. Hedetniemi, A.A. McRae and R. Mohan, Upper bounds on the coalition number, {\it Austral. J. Combin.} {\bf 80(3)} (2021),
	442--453.
	
	
	\bibitem{Abbas} A. Jafari, S. Alikhani, D. Bakhshesh, $k$-Coalitions in Graphs, 
	{\it Austral. J. Combin.}, {\bf 92(2)} (2025),  194-209.  
	
	
	\bibitem{Mojdeh1} 
	D.A. Mojdeh, I. Masoumi, Edge coalitions in graphs, Available at 
	\texttt{https://arxiv.org/abs/2302.10926}.
	
	\bibitem{Mojdeh2} 
	D.A. Mojdeh, M.R. Samadzadeh, Perfect coalition in graphs, Available at 
	\texttt{https://arxiv.org/pdf/2409.10185}.
	
	
	
	
\end{thebibliography}
\end{document}